\documentclass{amsart}

\input xy

\xyoption{all}

\usepackage{amssymb,amsbsy,amsthm,amsmath,graphicx,epsfig,times}

\newtheorem{thm}{Theorem}[section]
\newtheorem*{thm*}{Theorem}
\newtheorem{lem}[thm]{Lemma}
\newtheorem{prop}[thm]{Proposition}

\theoremstyle{remark}

\renewcommand{\rm}[1]{\mathrm{#1}}
\newcommand{\cal}[1]{\mathcal{#1}}

\newcommand{\bbR}{\mathbb{R}}

\newcommand{\bbZ}{\mathbb{Z}}

\newcommand{\sfP}{\mathsf{P}}

\renewcommand{\d}{\mathrm{d}}

\newcommand{\F}{\mathcal{F}}

\renewcommand{\O}{\Omega}

\renewcommand{\a}{\alpha}

\newcommand{\eps}{\varepsilon}

\renewcommand{\l}{\lambda}
\newcommand{\w}{\omega}

\newcommand{\fin}{\nolinebreak\hspace{\stretch{1}}$\lhd$}

\begin{document}

\title{A CAT($0$)-valued pointwise ergodic theorem}
\author{Tim Austin}
\date{}

\maketitle

\begin{abstract}
This note proves a version of the pointwise ergodic theorem for
functions taking values in a separable complete CAT($0$)-space.  The precise setting consists of an amenable locally compact group $G$ with left Haar
measure $m_G$, a jointly measurable, probability-preserving action $T:G\curvearrowright
(\O,\cal{F},\sfP)$ of $G$ on a probability space, and a separable complete CAT($0$)-space $(X,d)$ with barycentre map $b$.  In this setting we show that if $(F_n)_{n\geq 1}$ is a tempered F\o lner
sequence of compact subsets of $G$ and $f:\O\to
X$ is a measurable map such that for some (and hence any) fixed
$x \in X$ we have
\[\int_\O d(f(\w),x)^2\,\sfP(\d \w) < \infty,\]
then as $n\to\infty$ the functions of empirical barycentres
\[\w \mapsto b\Big(\frac{1}{m_G(F_n)}\int_{F_n}\delta_{f(T^g\w)}\,m_G(\d g)\Big)\]
converge pointwise for almost every $\w$ to a $T$-invariant function
$\bar{f}:\O\to X$.
\end{abstract}

\parskip 7pt

\section{Introduction}

Suppose that $(\O,\cal{F},\sfP)$ is a probability space and $(X,d)$
a complete separable CAT($0$)-space (see, for instance, Bridson and
Haefliger~\cite{BriHae99}). We write $L^2(\sfP;X)$ for the space of
all measurable maps $f:\O\to X$ such that for some fixed point $x
\in X$ we have
\[\int_\O d(f(\w),x)^2\,\sfP(\d\w) < \infty.\]
It is easy to see that in this case this actually holds for every $x
\in X$, and that if $f,g \in L^2(\sfP;X)$ then also
\[\int_\O d(f(\w),g(\w))^2\,\sfP(\d\w) < \infty.\]
If we now define
\[d_2(f,g) := \sqrt{\int_\O d(f(\w),g(\w))^2\,\sfP(\d\w)},\]
then this is a metric on $L^2(\sfP;X)$ that is easily seen to be
also complete and CAT($0$), and separable if $\F$ is countably
generated up to $\sfP$-negligible sets.

In addition, let $P_2(X)$ be the collection of probability measures
on $X$ with finite second moment, in that sense that $\mu \in
P_2(X)$ if for some (and hence every) $x\in X$ we have
\[\int_Xd(y,x)^2\,\mu(\d y) < \infty.\]
In these terms the condition that $f \in L^2(\sfP;X)$ is equivalent
to $f_\#\sfP \in P_2(X)$, where $f_\#\sfP$ is the pushforward
measure of $\sfP$ under $f$.

The important geometric property of complete CAT$(0)$-spaces that
motivates our work is that they support a sensible notion of
averaging. More specifically, it has been known essentially since
work of Cartan~\cite{Car51} that for any $\mu \in P_2(X)$ there is a
unique point $x \in X$ for which the above integral is minimized
(see Chapter II.2 of Bridson and Haefliger). We will refer to this
as the \textbf{barycentre} of $\mu$ and denoted it by $b(\mu)$.

In terms of these barycentres we can now define a CAT$(0)$-notion of
ergodic averages. Suppose that $G$ is an amenable locally compact
group with left-invariant Haar measure $m_G$ that acts on
$(\O,\F,\sfP)$ through a jointly measurable, $\sfP$-preserving
action $g\mapsto T^g$.  Given this, for a measurable subset $E
\subseteq G$ with $m_G(E) < \infty$ and a point $\w \in \O$ we will
write
\[\nu_{f,E}(\w)= \frac{1}{m_G(E)}\int_E \delta_{f(T^g\w)}\,m_G(\d g)\]
for the \textbf{empirical measure} of $f$ across the associated
orbit patch of $T$: more explicitly, this is defined by
\[\nu_{f,E}(\w)(A) := \frac{1}{m_G(E)}m_G\{g\in E:\ f(T^g\w)\in A\}.\]
It is natural to view the function of barycentres of the empirical
measures of $f$,
\[\w \mapsto b(\nu_{f,E}(\w)),\]
as a CAT($0$) analog of the ergodic averages
\[\frac{1}{m_G(E)}\int_Ef(T^g\w)\,m_G(\d g)\]
available in case $f:\O\to\bbR$.

In the case of real-valued functions, it is known that if $f \in
L^2(\sfP)$ then these ergodic averages converge for $\sfP$-almost
every $\w \in \O$ as $E$ increases along a suitably chosen F\o lner
sequence of subsets of $G$.  This follows from the classical
pointwise ergodic theorem of Birkhoff in case $G = \bbZ$, and has
also long been known for many other concrete groups such as $\bbZ^d$
or $\bbR^d$. The general case was established only quite recently by
Lindenstrauss, who found that the appropriate condition to place on
the F\o lner sequence $(F_n)_{n\geq 1}$ is that it be
\textbf{tempered}: this holds if for some fixed $C > 0$ and all
$n\geq 1$ we have
\[m_G\Big(\bigcup_{k < n}F_k^{-1}F_n\Big) \leq Cm_G(F_n)\]
(this is also referred to as the `Shulman condition').

By thinning out an initially-given F\o lner sequence if necessary it
follows that any amenable group does admit F\o lner sequences
satisfying this condition.  Given such a sequence, Lindenstrauss
proves (alongside other results) that for any $f\in L^1(\sfP)$ there
is a $T$-invariant function $\bar{f}:\O\to \bbR$ such that
\[\frac{1}{m_G(F_n)}\int_{F_n}f(T^g\w)\,m_G(\d g) \to \bar{f}(\w)\]
for $\sfP$-almost every $\w \in \O$.

In the present note we will show that this result can be extended to
maps in $L^2(\sfP;X)$ for a separable CAT$(0)$-space $X$, replacing
ergodic averages with the empirical measure barycentres introduced
above:

\begin{thm}\label{thm:main}
If $T:G \curvearrowright (\O,\F,\sfP)$, $(X,d)$ are as above, then
for any $f \in L^2(\sfP;X)$ there is a $T$-invariant function
$\bar{f}:\O\to X$ such that for any tempered F\o lner sequence
$(F_n)_{n\geq 1}$ of compact subsets of $G$ we have
\[b(\nu_{f,E}(\w))\to \bar{f}(\w)\]
for $\sfP$-a.e. $\w \in \O$.

In particular, if $T$ is ergodic then $\bar{f}$ is just the constant
function with value $b(f_\#\sfP) \in X$.
\end{thm}

\textbf{Remark}\quad We restrict to separable $X$ in order to avoid
discussing the nuances between different notions of `measurability'
for $f$, but provided the right notion is chosen this seems to make
no real restriction. \fin

We will find that this theorem follows quite quickly from an appeal
to the real-valued pointwise ergodic theorem, together with an
approximation argument based on a maximal ergodic theorem for such
group actions and F\o lner sequences also obtained by Lindenstrauss
in~\cite{Lin01}.  On the other hand, since any Hilbert space is
CAT$(0)$ with barycentre map simply given by averaging,
Theorem~\ref{thm:main} does contain the pointwise ergodic theorem
for square-integrable maps into a separable Hilbert space as a
special case.

It is more subtle to find a theorem about CAT$(0)$ targets that encompasses the real-valued theorem for arbitrary functions in $L^1(\sfP)$, because the condition that $\mu
\in P_2(X)$ already appears in the definition of the barycentre
$b(\mu)$.  In fact, suitable notions of barycentre defined for all measures with finite first moment have been set up by Es-Sahib and Heinich~\cite{EsSHei99} and by Sturm~\cite{Stu03}, but they require more delicate handling than the map $b$ we use here. Since the present note was submitted, I understand that Bruno Duchesne and Andr\'es Navas have independently shown that the properties of these alternative barycentre maps do lead to a CAT$(0)$-valued $L^1$-ergodic theorem~(\cite{Dus10,Nav11}), essentially via the argument below.

In the present generality, Theorem~\ref{thm:main} is new even in the classical case $G
= \bbZ$, but it does have a precedent for i.i.d. sequences of $X$-valued random variables with finite second moment.  For such sequences the pointwise convergence of empirical measures (that is, the Law of Large Numbers) was shown by Sturm in~\cite{Stu03}: see his Proposition 6.6, and compare also with his Theorem 4.7, in which a related convergence result is used to give an alternative characterization of the barycentre map.  The case of i.i.d. sequences of random variables corresponds to an action $\bbZ\curvearrowright (\O,\cal{F},\sfP)$ and function $f:\O\to X$ for which all the distinct composites $f\circ T^n$ are independent.  It would be interesting to know whether the greater generality of Theorem~\ref{thm:main} yields any additional applications in the study of group actions on CAT$(0)$ spaces or related structures, for example through generalizing these results of Sturm, or in connexion with recent results identifying barycentres as maximum likelihood estimators in statistics (see, for instance, Dewarrat and Ruh~\cite{DewRuh02} and Fl\"uge and Ruh~\cite{FluRuh06}).

Theorem~\ref{thm:main} also bears comparison with various other works
relating probability-preserving actions to the geometry of
CAT$(0)$-spaces; we refer the reader in particular to the proof by
Karlsson and Margulis in~\cite{KarMar99} of an analog of Oseledets'
Theorem for cocycles over a probability-preserving action taking
values in the semigroup of contractions of a CAT$(0)$-space.

In the next section we will recall the maximal ergodic theorem we
need from Lindenstrauss~\cite{Lin01} and derive from it a useful
consequence for CAT$(0)$-valued maps, and then in
Section~\ref{sec:mainproof} we will use these results to prove
Theorem~\ref{thm:main}.

\textbf{Acknowledgements}\quad My thanks go to Lior Silberman for
helpful suggestions, to Andr\'es Navas for pointing out an error in an earlier version, and to the anonymous referees for some relevant references. \fin

\section{A CAT$(0)$-valued maximal ergodic theorem}\label{sec:max}

An important innovation behind Lindenstrauss' proof of his pointwise
ergodic theorem was a proof of the weak-$(1,1)$ maximal ergodic
theorem in the setting of general amenable groups and tempered F\o
lner sequences.

Let us introduce the standard notation $Mf$ for the ergodic maximal
function associated to a F\o lner sequence and a function $f:\O\to
\bbR$:
\[Mf(\w) := \sup_{n\geq 1}\frac{1}{m_G(F_n)}\int_{F_n}f(T^g\w)\,m_G(\d g).\]
(see, for example, Peterson's book~\cite{Pet83} for background on
this maximal function in the classical case $G = \bbZ$, or Section 3
of Lindenstrauss~\cite{Lin01}).

\begin{prop}[Theorem 3.2 of~\cite{Lin01}]\label{prop:max-erg-th}
Let $T:G\curvearrowright (\O,\F,\sfP)$ be an action as above and
$(F_n)_{n\geq 1}$ a tempered F\o lner sequence.  Then there is a $c
> 0$, depending on the sequence $(F_n)_{n\geq 1}$ but not on $X$,
such that for any $f \in L^1(\sfP)$ we have
\[\sfP\{\w\in \O:\ Mf(\w) > \a\} \leq \frac{c}{\a}\|f\|_1\]
for all $\a > 0$. \qed
\end{prop}

We will make use of this result to prove the following analog for
CAT$(0)$-valued maps.

\begin{thm}[CAT$(0)$-valued maximal ergodic
theorem]\label{thm:CAT0-max-erg-th} If $f,h:\O\to X$ are two members
of $L^2(\sfP;X)$ then
\[\sfP\big\{\w\in\O:\ \sup_{n\geq 1}d(b(\nu_{f,F_n}(\w)),b(\nu_{h,F_n}(\w))) > \a\big\} \leq \frac{c}{\a^2}d_2(f,h)^2\]
for every $\a \in (0,\infty)$, where $c$ is the same constant as in
Proposition~\ref{prop:max-erg-th}.
\end{thm}

In order to prove this, we need an elementary result controlling the
behaviour of barycentres in terms of the Wasserstein metric on
$P_2(X)$.  Recall that this is defined for $\mu,\nu \in P_2(X)$ by
\[W_2(\mu,\nu) := \inf_{\l\ \rm{a}\ \rm{coupling}\ \rm{of}\ \mu,\nu}\sqrt{\int_{X^2}d(x,y)^2\,\l(\d x,\d y)}.\]
As is standard, this defines a metric on $P_2(X)$ (see, for
instance, Villani~\cite{Vil03}).  In the presence of the
CAT$(0)$ condition, it controls barycentres as follows.

\begin{lem}\label{lem:Wasser}
If $\mu,\nu \in P_2(X)$ then
\[d(b(\mu),b(\nu))\leq W_2(\mu,\nu).\]
\end{lem}

\textbf{Remark}\quad In the published version of this paper, the proof given for Lemma~\ref{lem:Wasser} was incorrect.  A corrected proof has been published as an erratum in the same journal.  This preprint has been re-written with the correct proof.  It actually gives improved control by the weaker Wasserstein metric $W_1$, rather than $W_2$, but I have not changed the remainder of the preprint to account for this. I am grateful to David Fisher for bringing the mistake in the original paper to my attention, and to Assaf Naor for pointing me to the relevant reference~\cite{LanPavSch00}. \fin

%

\textbf{Proof}\quad 
This can be deduced from a small generalization of a known result:~\cite[Lemma 4.2]{LanPavSch00}, which gives control of a similar kind for measures supported on finitely many atoms.  We quickly repeat the argument from~\cite{LanPavSch00} in the generality we need.

For any $x,y \in X$, we have
\[d(x,b(\nu))^2 + d(y,b(\mu))^2 \leq d(x,b(\mu))^2 + d(y,b(\nu))^2 + 2d(b(\mu),b(\nu))d(x,y)\]
(this is a standard inequality for any four points in a CAT$(0)$-space: see~\cite[Lemma 2.1]{LanPavSch00}). Letting $\l$ be any coupling of $\mu$ and $\nu$, we may integrate this inequality with respect to $\l(\d x,\d y)$ to obtain
\begin{multline}\label{eq:ineq}
\int d(x,b(\nu))^2\,\mu(\d x) + \int d(y,b(\mu))^2\,\nu(\d y)\\
\leq \int d(x,b(\mu))^2\,\mu(\d x) + \int d(y,b(\nu))^2\,\nu(\d y) + 2d(b(\mu),b(\nu))\int d(x,y)\,\l(\d x,\d y).
\end{multline}
On the other hand, another standard inequality for baycentres in CAT$(0)$-spaces (see~\cite[Lemma 4.1]{LanPavSch00}) gives that
\[\int d(x,b(\nu))^2\,\mu(\d x) \geq \int d(x,b(\mu))^2\,\mu(\d x) + d(b(\mu),d(\nu))^2,\]
and similarly with the roles of $\mu$ and $\nu$ reversed.  Inserting these two inequalities into~(\ref{eq:ineq}), we obtain
\begin{multline*}
\int d(x,b(\mu))^2\,\mu(\d x) + \int d(y,b(\nu))^2\,\nu(\d y) + 2d(b(\mu),b(\nu))^2\\
\leq \int d(x,b(\mu))^2\,\mu(\d x) + \int d(y,b(\nu))^2\,\nu(\d y) + 2d(b(\mu),b(\nu))\int d(x,y)\,\l(\d x,\d y).
\end{multline*}
Upon simplifying, this yields
\[d(b(\mu),b(\nu)) \leq \int d(x,y)\,\l(\d x,\d y) \leq \sqrt{\int d(x,y)^2\,\l(\d x,\d y)},\]
where the second bound follows by the Cauchy--Bunyakowski--Schwartz inequality.  Since $\l$ was an arbitrary coupling, this completes the proof. \qed

\textbf{Proof of Theorem~\ref{thm:CAT0-max-erg-th}}\quad For any $\a
> 0$ the above lemma implies that
\begin{multline*}
\{\w \in \O:\ \sup_{n\geq 1}d(b(\nu_{f,F_n}(\w)),b(\nu_{h,F_n}(\w)))
> \a\}\\ \subseteq \{\w \in \O:\ \sup_{n\geq
1}W_2(\nu_{f,F_n}(\w),\nu_{h,F_n}(\w)) > \a\}.
\end{multline*}
On the other hand, for each $n$ the measure
\[\l_n := \frac{1}{m_G(F_n)}\int_{F_n}\delta_{(f(T^g\w),h(T^g\w))}\,m_G(\d g)\]
on $X^2$ is clearly a joining of $\nu_{f,F_n}(\w))$ and
$\nu_{h,F_n}(\w)$, and so
\[W_2(\nu_{f,F_n}(\w),\nu_{h,F_n}(\w))^2 \leq \frac{1}{m_G(F_n)}\int_{F_n}d(f(T^g\w),h(T^g\w))^2\,m_G(\d g).\]

If we now write $F(\w) := d(f(T^g\w),h(T^g\w))^2$, then combining
the above observations gives
\[\{\w \in \O:\ \sup_{n\geq
1}d(b(\nu_{f,F_n}(\w)),b(\nu_{h,F_n}(\w))) > \a^2\} \subseteq \{\w
\in \O:\ MF(\w) > \a\},\] and so now
Proposition~\ref{prop:max-erg-th} gives
\[\sfP\{\w \in \O:\ \sup_{n\geq
1}d(b(\nu_{f,F_n}(\w)),b(\nu_{h,F_n}(\w))) > \a\} \leq
\frac{c}{\a^2}\|F\|_1\] and so finally observing that $\|F\|_1 =
d_2(f,h)^2$ completes the proof. \qed

\section{Proof of the main theorem}\label{sec:mainproof}

\textbf{Proof of the Theorem~\ref{thm:main}}\quad Let us first prove
the assertion of Theorem~\ref{thm:main} for a finite-valued function
$h:\O\to X$.  To this is associated some finite measurable partition
$\O = A_1 \cup A_2 \cup \ldots \cup A_m$ and collection of points
$x_1,x_2,\ldots,x_k \in X$ such that $h(\w) = x_i$ when $\w \in
A_i$. This case now follows easily from Lemma~\ref{lem:Wasser} and
the real-valued pointwise ergodic theorem: given $\eps > 0$, for
almost every $\w \in \O$ and every $i \leq k$, that theorem gives
some $n(\w,k,\eps) \geq 1$ such that
\[n\geq n(\w,k,\eps) \quad\quad\Rightarrow\quad\quad
\Big|\frac{1}{m_G(F_n)}m_G\{g\in F_n:\ T^g\w\in A_i\} -
\sfP_\omega(A_i)\Big| < \eps\] for some $T$-invariant function $\omega \mapsto \sfP_\omega(A_i)$ (where, of course, $\sfP_\omega$ is actually just the component over $\omega$ of the decomposition of $\sfP$ into $T$-ergodic components). Choosing $n$ sufficiently large we
can make this hold for this $\eps$ and $\w$ and all $i \leq k$. It
follows that the empirical measure
\[\nu_{h,F_n}(\w) = \frac{1}{m_G(F_n)}\int_{F_n}\delta_{h(T^g\w)}\,m_G(\d g)\]
satisfies the total variation inequality
\[\|\nu_{h,F_n}(\w) - \sfP_\omega(A_1)\delta_{x_1} - \sfP_\omega(A_2)\delta_{x_2} - \cdots - \sfP_\omega(A_k)\delta_{x_k}\|_{\rm{var}} < \eps,\]
and hence that
\[W_2\big(\nu_{h,F_n}(\w),\,\sfP_\omega(A_1)\delta_{x_1} + \sfP_\omega(A_2)\delta_{x_2} + \cdots + \sfP_\omega(A_k)\delta_{x_k}\big) < \sqrt{\eps}\,\rm{diam}\{x_1,x_2,\ldots,x_k\}.\]

Combining this with Lemma~\ref{lem:Wasser} we obtain that for
almost every $\w$, for every $\eps
> 0$ we have
\begin{multline*}
d(b(\nu_{h,F_n}(\w)),b(\sfP_\omega(A_1)\delta_{x_1} + \sfP_\omega(A_2)\delta_{x_2} + \cdots + \sfP_\omega(A_k)\delta_{x_k}))\\ < \sqrt{\eps}\,\rm{diam}\{x_1,x_2,\ldots,x_k\}
\end{multline*}
for all sufficiently large $n$, and so since
$\rm{diam}\{x_1,x_2,\ldots,x_k\}$ is a fixed quantity for a given
$h$ it follows that
\[b(\nu_{h,F_n}(\w))\to b(\sfP_\omega(A_1)\delta_{x_1} + \sfP_\omega(A_2)\delta_{x_2} + \cdots + \sfP_\omega(A_k)\delta_{x_k})\]
as $n\to\infty$ for almost every $\w$.

Finally, the same assertion for an arbitrary $f \in L^2(\sfP;X)$
follows from a routine approximation argument and appeal to
Theorem~\ref{thm:CAT0-max-erg-th}.  Letting $c
> 0$ be the constant of Theorem~\ref{prop:max-erg-th} and given $\a
> 0$, since $X$ is separable we can always find a finite-valued function $\psi:\O\to
X$ such that $d_2(f,h) < \a^2$, and hence
\[\sfP\big\{\w\in\O:\ \sup_{n\geq 1}d(b(\nu_{f,F_n}(\w)),b(\nu_{h,F_n}(\w))) > \a\big\} \leq \frac{c}{\a^2}d_2(f,h)^2 < c\a^2,\]
so since we have seen that $b(\nu_{h,F_n}(\w))$ tends to a limit that is a $T$-invariant function of $\omega$ on a
conegligible set of $\w \in \O$ it follows that outside the above
subset of $\O$ of measure at most $c\a^2$, the sequence
$(b(\nu_{f,F_n}(\w)))_{n=1}^\infty$ asymptotically oscillates by at
most $\a$ in $X$, and so since $\a$ was arbitrary this sequence must
actually converge for almost every $\w$ to a function which is also $T$-invariant.  This completes the proof.
\qed

\textbf{Remark}\quad From another simple approximation by a
finite-valued function as in the above proof and an appeal to the
Dominated Convergence Theorem, it follows directly from
Theorem~\ref{thm:main} that the empirical barycentres of $f$ also
converge to $\bar{f}$ in the metric space $(L^2(\sfP;X),d_2)$. \fin

\bibliographystyle{abbrv}
\bibliography{bibfile}

\parskip 0pt

\vspace{7pt}

\small{\textsc{Department of Mathematics, Brown University,
Providence, RI 02912, USA}

Email: \verb|timaustin@math.brown.edu|

URL: \verb|http://www.math.brown.edu/~timaustin|}

\end{document}